\documentclass[letterpaper,11pt]{article}
\usepackage{times}
\usepackage{amsmath, amsthm, amssymb}    
\usepackage{bridges}			
\usepackage{graphicx}			
\usepackage{wrapfig}
\usepackage{subcaption}			

\usepackage[utf8]{inputenc}
\usepackage[T1]{fontenc}

\usepackage[colorlinks=true, urlcolor=blue, citecolor=black, linkcolor=black]{hyperref}  
\graphicspath{{figures/}}

\title{The Path to Aperiodic Monotiles\footnote{I was invited to submit
an article to \textit{Eureka}, the journal of the Archimedeans (the student
mathematics society at Cambridge University).  This article appeared in
Issue 67 of the journal; I have made one small update, to reflect a new
lower bound on the undecidability of the tiling problem~\cite{DL2025}.  
For more information about \textit{Eureka}, visit \url{https://archim.soc.srcf.net/?page_id=140}.}}

\author{Craig S. Kaplan \\ \\
School of Computer Science, University of Waterloo, Ontario, Canada; csk@uwaterloo.ca} 

\date{}					

\begin{document}

\maketitle

In the 1978 issue of this journal (and in a reprint in 2012), Sir 
Roger Penrose told the story of the sets of aperiodic tiles that now
bear his name.  These sets, the most famous of which consists of
two shapes called the ``kite'' and ``dart'', have the remarkable
property that they tile the plane, but never periodically.  


Penrose recognized that his kite and dart represented a leap forward
in the search for small aperiodic tile sets, but that one final
step was theoretically possible: at the end of his essay he observed
that ``it is not known whether there is a \emph{single} shape that
can tile the Euclidean plane non-periodically.''  The earliest
instance of this question that I can find in print comes from Martin
Gardner's Mathematical Games column in January 1977, which first
popularized Penrose tiles.  When the column was reprinted in 1989
in Gardner's anthology \emph{From Penrose Tiles to Trapdoor Ciphers},
he expanded the text to refer to this question as ``the major
unsolved problem'' in tiling theory, adding that most mathematicians
doubted such a shape could exist.

In late 2022, David Smith, an online acquaintance from the world
of tiling enthusiasts, emailed me a drawing that would change our
understanding of aperiodicity, not to mention my life.  David's
shape, which we now call the ``hat'', is a positive answer to the
problem posted by Penrose and Gardner.  It is an \emph{aperiodic
monotile}: a single shape that can only tile the plane non-periodically.
David and I collaborated with Chaim Goodman-Strauss and Joseph
Samuel Myers (a Cambridge alumnus and Archimedean) to publish a
proof of the hat's aperiodicity~\cite{hat}.

It is a privilege to be offered a place in these pages to relate
one more chapter in the story of tiling theory.  I will arrive at
a discussion of aperiodic monotiles along a route that visits two
favourite topics of mine that are closely related: Heesch numbers
and isohedral numbers.  I and others always viewed the study of
these two topics as a means of making incremental progress towards
aperiodicity.  Despite David's discovery of the hat, Heesch numbers
and isohedral numbers still offer a lot of independent interest
(and in hindsight can now perhaps be regarded as even more difficult
problems).  Throughout, I will refer liberally to unsolved problems,
both established and new, which I hope will induce some readers to
delve more deeply into the study of tiling theory.

%

I require just a few basic mathematical notions to establish context
for the discussion that follows.  A \emph{shape} is a topological
disk.  Given a set of shapes $\mathcal{S}=\{S_1,\ldots,S_k\}$, a
\emph{tiling} from $\mathcal{S}$ is a countably infinite collection
of tiles $\mathcal{T} = \{T_1,T_2,\ldots\}$, where each $T_i$ is
congruent to a shape in $\mathcal{S}$, and where the tiles cover
the entire plane with no gaps and no overlaps (except on their
boundaries).  We also say that $\mathcal{S}$ \emph{admits} the
tiling~$\mathcal{T}$.  A \emph{monotile} is a shape that admits a
tiling as a singleton set.  A \emph{patch} is a finite set of
non-overlapping shapes whose union is a topological disk.

Two shapes are \emph{neighbours} if they have at least one point
in common on their boundaries.  We frequently make a simplifying
``finite neighbour'' assumption: given shapes $S_1$ and $S_2$, a
congruent copy of $S_1$ can be neighbours with a congruent copy of
$S_2$ in only finitely many different ways (or, in the context of
monotiles, that two copies of a shape $S$ can be neighbours in only
finitely many different ways).  This assumption tends to ``discretize''
the behaviour of shapes in tilings and patches, making them more
amenable to computation.  We can usually remove this assumption
later, sometimes altering the outlines of shapes to enforce the
neighbour relationships we want.  

\begin{figure}
\begin{center}
\includegraphics[width=5in]{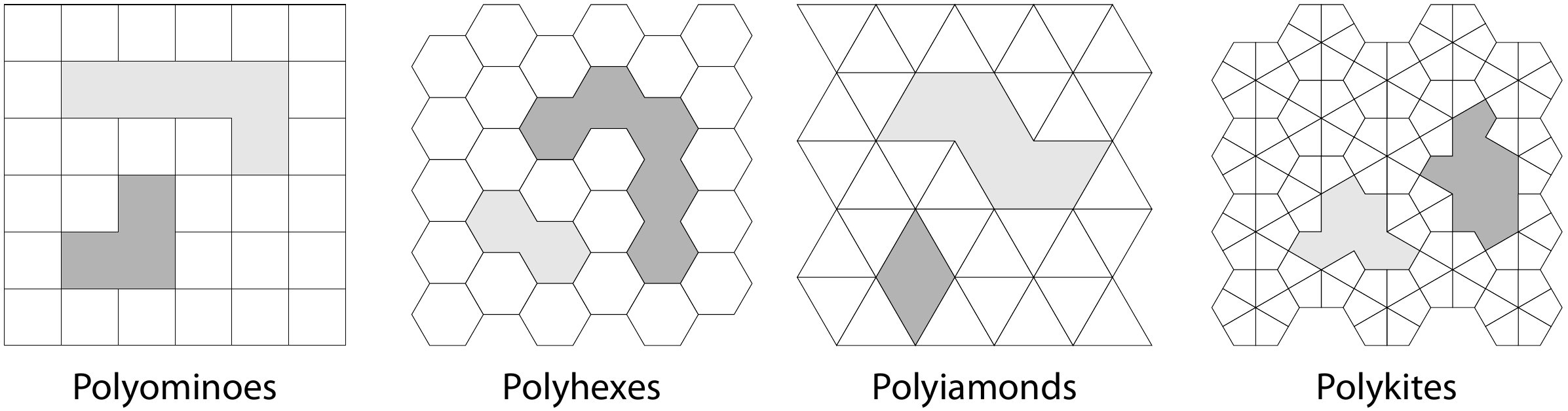}
\end{center}
\caption{\label{fig:polyforms}Examples of polyforms.}
\end{figure}

Specifically, we often restrict attention to \emph{polyforms},
shapes constructed by gluing together copies of some basic unit
cell from an underlying tiling (Figure~\ref{fig:polyforms}).  For
example, \emph{polyominoes}, \emph{polyhexes}, and \emph{polyiamonds}
are glued-together squares, regular hexagons, and equilateral
triangles, respectively.  The hat is a \emph{polykite}: its units
are kite-shaped cells from a tiling created by overlaying a regular
hexagon tiling and its dual equilateral triangle tiling.  A polyform
naturally has finite neighbours if we assume that copies of it  must
be placed in alignment with the cell tiling from which it was
constructed.  Polyforms offer a convenient means of probing broader
questions about tilings.  Unlike general shapes they can be enumerated,
and many of their tiling theoretic properties can be calculated
algorithmically.  In focusing upon polyforms, we make the optimistic
(and thus far fruitful) assumption that they embody many of the
interesting properties we would like to understand about tilings
in general.

\section*{Heesch Numbers}

Given a shape with unknown tiling theoretic properties, the most
obvious first question to ask is whether it admits any tilings at
all.  It is far from clear how one might go about answering that
question in a generic way.  Consider that a given shape might admit
no tilings, a finite number of distinct tilings, or uncountably
many distinct tilings (even when restricted to finite neighbours).
The problem of determining whether a shape tiles the plane may even
be computationally undecidable, in the sense that a Turing machine
can be converted into a shape that tiles the plane if and only if
the Turing machine runs forever.  The tiling problem is known to
be undecidable for as few as three polygons~\cite{DL2025},
but the question remains
open for smaller sets.  If the problem were undecidable for a single
shape, then roughly speaking the question of whether a shape tiles
the plane would be ``as hard as possible'', and unanswerable in general
terms.

Absent a comprehensive algorithm, we might at least devise some simple
tests to determine a shape's behaviour.  For example, we could try
to surround a shape by copies of itself; if a shape cannot even be
surrounded, then surely it does not tile the plane.  Specifically,
we say that a set of shapes $\{S_1,\ldots,S_k\}$ \emph{surrounds}
a shape $S$ if all these shapes together form a patch with $S$ in its
interior, such that every $S_i$ is a neighbour of $S$.  (Here I ignore
the possibility that every part of the boundary of $S$ is covered by
one of the~$S_i$, but that the union of the shapes contains a hole that
can be filled by additional copies of $S$.)

If a shape $S$ can meet copies of itself in only finitely many ways,
then an algorithm for checking its surroundability clearly exists:
a surround must be a subset of the finite set of neighbours of $S$,
and all such subsets can in principle be checked to see if any is
a surround.  This algorithm is obviously inefficient, but I am not
aware of an alternative test of surroundability whose computational
complexity is asymptotically better.  While a tractable algorithm
would be convenient, I believe that a more natural target for new
research would be to prove that surroundability is NP-complete,
perhaps focusing on restricted classes of shapes like polyominoes.

While an unsurroundable shape clearly does not admit tilings, the
converse does not hold.  As Heinrich Heesch first observed in 1968,
there exist shapes that can be surrounded once but not twice.
Let us introduce more terminology before inquiring how far
this process of surrounding may be extended.  Given a shape $S$, a
\emph{$0$-patch} is a single copy of $S$.  Now recursively define
a \emph{$(k+1)$-patch} to be a $k$-patch that is itself surrounded
by copies of $S$.  The \emph{Heesch number} of $S$ is then the
largest integer $n$ for which $S$ has an $n$-patch.  If $S$ admits
a tiling, its Heesch number is defined to $\infty$.  A deep result
in tiling theory guarantees that if a shape does not tile the plane,
then its Heesch number must be finite (ruling out, for example, a
shape that tiles a quadrant but not the whole plane).

If there were a maximum finite Heesch number, then we could check in
finite time whether a shape $S$ admits a tiling.  Given $S$, we
recursively attempt to construct all possible $k$-patches for
increasing $k$.  If we ever reach a $k$ for which no $k$-patch
can be surrounded to create a $(k+1)$-patch, then $S$ does not admit
tilings.  Or, if we exceed the supposed maximum Heesch number,
then $S$ must admit a tiling, even if we know nothing about the
structure of that tiling.  Of course, if Heesch numbers are unbounded
then this approach is of limited use, because there is no way to know
when to abandon the search and assume that a shape admits tilings.

\begin{figure}
\begin{center}
\includegraphics[width=\textwidth]{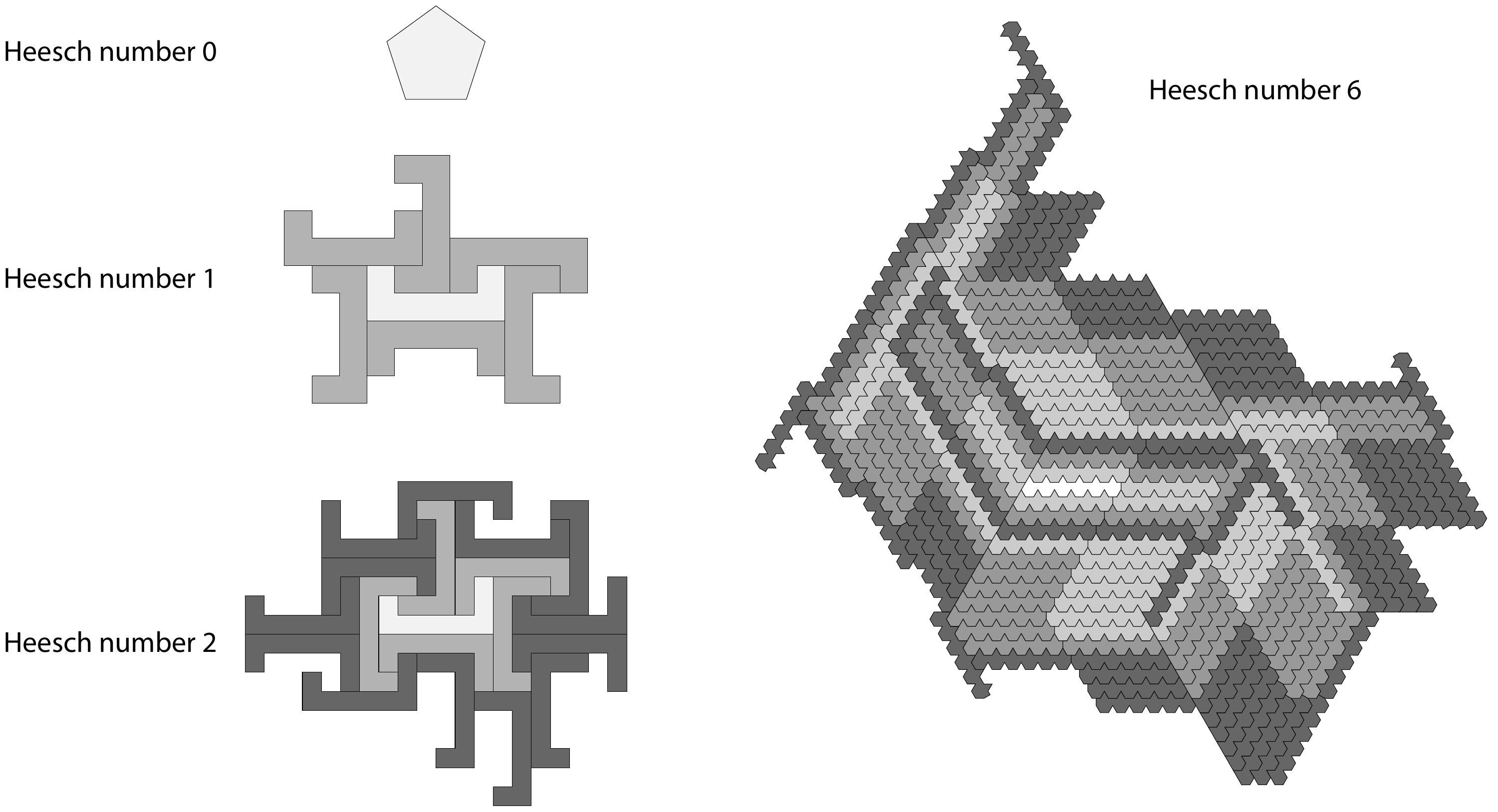}
\end{center}
\caption{\label{fig:heesch}A few non-tiling shapes with their 
	Heesch numbers.}
\end{figure}

\emph{Heesch's problem} asks which positive integers are Heesch
numbers, or more broadly whether they are bounded.  Very little is
known about limitations on Heesch numbers.  We have examples of
shapes with finite Heesch numbers from zero up to a rather paltry
six (Figure~\ref{fig:heesch}),
with the remarkable record holder published in 2021 by
Bojan Ba\v{s}i\'c~\cite{Basic2021}.  

Heesch's problem is one of my favourite unsolved
problems in mathematics.  Other problems like the Riemann Hypothesis
and P versus NP are more consequential, but in both cases everyone
knows what the ``right'' answer is, even if the proofs remain
elusive.  With Heesch's problem, I am fully prepared to be astounded
by any outcome.  I am inclined to believe that Heesch numbers
have no bound. Even so, I cannot imagine how one might define, say,
a sequence of shapes $\{S_n\}_{n=0}^\infty$ such that each $S_n$
has Heesch number~$n$.  But it is even harder to conceive of a
reason why the Euclidean plane would impose an upper limit on Heesch
numbers of six, or 89, or 1000 (and my degree of astonishment would
grow in proportion to that limit, whatever it turned out to be).

In previous work, I wrote software to compute Heesch numbers of
non-tiling polyominoes, polyhexes, and polyiamonds, with the as-yet
unrealized goal of finding new record breakers~\cite{Kaplan}.  I
avoided the pathologically inefficient recursive construction
suggested above.  Instead, I expressed the question of whether a
polyform $S$ has an $n$-patch as a formula in propositional logic,
and used a type of software library known as a \emph{SAT solver}
to determine whether the formula could be made true through a
suitable choice of values for its variables.  Any satisfying
assignment could be translated back into a description of an
$n$-patch.  This approach offers no performance guarantees, but SAT
solvers employ powerful heuristics that work well in practice.  My
exhaustive calculations support the observation that the proportion
of shapes with Heesch number $n$ falls off precipitously as $n$
increases.  If nearly all non-tilers to be tested will have Heesch
number zero, then new research should seek full algorithms or
heuristics that reveal unsurroundability as quickly as possible,
avoiding the SAT solver in most cases. We can afford the larger
cost of the SAT solver in the small proportion of shapes where it
is needed.

\section*{Isohedral Numbers}

Checking that a shape is unsurroundable, i.e., that it has a Heesch
number of zero, can quickly identify many non-tilers.  We can further
narrow our inconclusive results by checking whether a shape's Heesch
number is less then some arbitrary threshold, and giving up if we
pass that threshold.  The flipside to such an approach would be to
peform some simple tests that establish definitively that a shape
\emph{does} tile. While we cannot account for all possible tilings
that a shape might admit, some types of tilings are simple enough and
common enough that they can be tested for explicitly.

Suppose that a monotile $S$ admits a tiling
$\mathcal{T}=\{T_1,T_2,\ldots\}$.  As a drawing in the plane, this
tiling has a symmetry group, consisting of rigid motions of the
plane that map the tiling to itself.  Given any two tiles $T_i$ and
$T_j$ in $\mathcal{T}$, we may ask whether one of those symmetries
happens to map $T_i$ to $T_j$.  If so, we refer to those two tiles
as \emph{transitively equivalent}. This equivalence relation
partitions the tiles in a tiling into transitivity classes, consisting
of the orbits of the tiles under the action of the tiling's symmetry
group.  A tiling is periodic, meaning that its symmetry group is
one of the 17 wallpaper groups that have translation symmetries in
two directions, if and only if its tiles fall into  a finite number
of transitivity classes.  A tiling is \emph{isohedral} if all its
tiles belong to a single transitivity class.  In an isohedral tiling,
the tiles are indistinguishable: every tile offers the same view
out to infinity in every direction.

Isohedral tilings are highly organized, so much so that we can go
beyond merely checking whether a shape admits an isohedral tiling:
we can often do so with provable efficiency.  In an isohedral tiling,
every tile must be surrounded by its immediate neighbours in the
same way, and that surround determines the complete structure of
the tiling.  It is therefore possible to devise algorithms that
check for isohedral tilability using purely local information about
a shape and its immediate neighbours.  The state of the art for
such algorithms is one by Langerman and Winslow that checks whether
a polyomino with perimeter $n$ tiles isohedrally in 
$\mathrm{O}(n\log^2n)$ time~\cite{LW2016}.

But as with Heesch numbers, checking whether a shape tiles isohedrally
cannot be the whole story: there exist shapes that tile the plane
periodically but never isohedrally.  The first such polygon was
exhibited by Heesch in 1935.  Accordingly, we introduce a number
that characterizes how disorderly a shape is forced to be in terms
of the periodic tilings it admits.  Define a shape's \emph{isohedral
number} to be the minimum number of transitivity classes in any
tiling by that shape.  To a first approximation, the isohedral
number measures how many copies of a shape you must glue together
before you obtain a patch that tiles in a simple way (i.e.,
isohedrally).  In practice the isohedral number is usually equal
to the size of smallest patch that tiles isohedrally, but the numbers
can diverge slightly when the shape and its patch are symmetric.
At a high level, every shape that admits a periodic tiling has a
finite isohedral number.

Isohedral numbers are the cousins of Heesch numbers, from the
opposite shore of the tiling theoretic sea.  As with Heesch numbers,
we do not know which positive integers can be isohedral numbers,
or whether these numbers are bounded.  If there were a bound on
isohedral numbers, then an algorithm would exist to determine whether
a shape admits a periodic tiling: we could construct all patches
up to a predetermined size and check whether any of them tiles
isohedrally.

\begin{figure}
\begin{center}
\includegraphics[width=\textwidth]{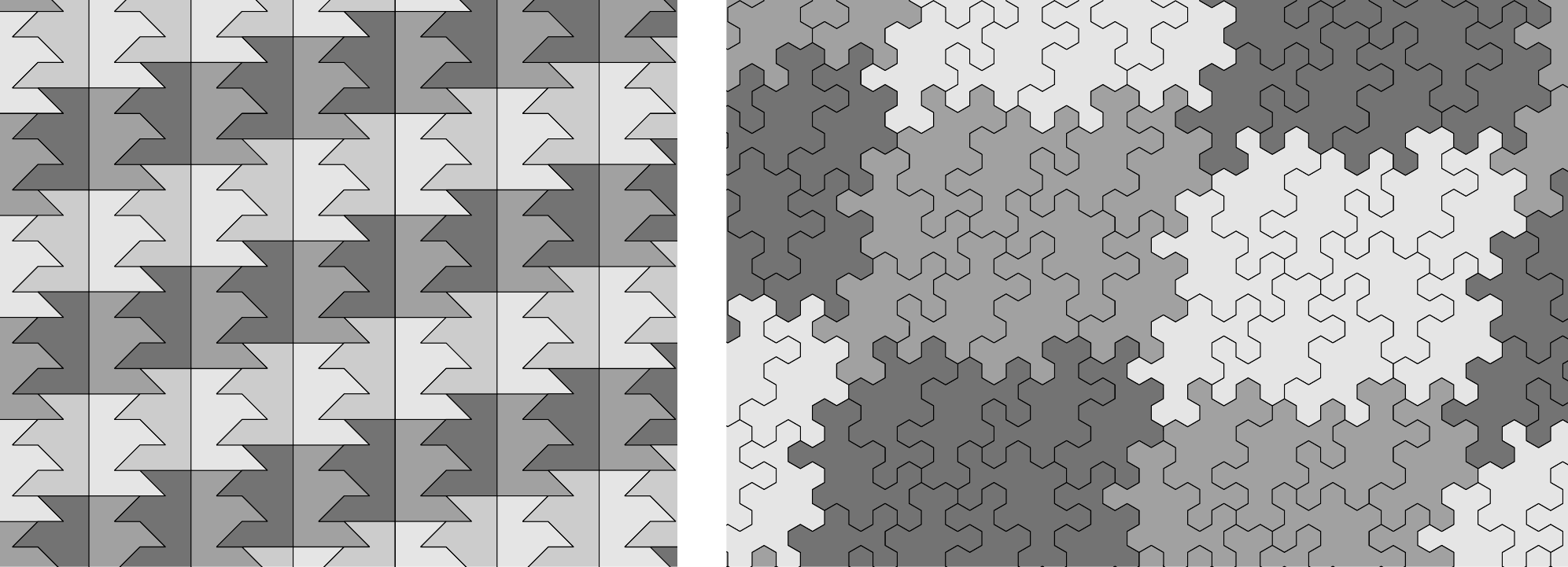}
\end{center}
\caption{\label{fig:anisohedral}Heesch's tile with isohedral
	number 2 (left) and Myers's 16-hex with isohedral number 10 (right).}
\end{figure}

Joseph Myers is the main source of empirical data on this subject.
He conducted an extensive computation of isohedral numbers of
polyominoes, polyhexes, polyiamonds, and polykites~\cite{Myers}.
He discovered polyforms with every isohedral number up to ten
(Figure~\ref{fig:anisohedral}), excluding seven.  His software
ultimately relies on brute-force search to construct patches that
tile isohedrally.  However, high isohedral numbers are so rare, and
Joseph's code is so well optimized for the common cases of low
isohedral numbers and non-tilers, that his software is lightning
fast on average.

The lack of an example in Joseph's data of a shape with isohedral
number seven is a hole begging to be filled.  More generally,
we must venture onward in search of shapes with ever higher isohedral
numbers, or articulate some mathematical reason why they must be
bounded.  I see no reason why a bound should exist; if one does, I
would be extremely surprised if it were ten.

\section*{Aperiodicity}

Armed with algorithms for computing Heesch numbers and isohedral
numbers, we might imagine combining them into a master algorithm
that climbs those two ladders in tandem.  Given a shape $S$
we check, for each successive integer $n$, whether
$S$ has Heesch number at most $n-1$ (by failing to construct
an $n$-patch) or isohedral number at most $n$ (by constructing a
patch of $n$ shapes that tiles isohedrally).  As soon as either
of these tests succeeds, we stop and report that $S$ is a non-tiler
with a given Heesch number or a periodic tiler with a given isohedral
number.

In order to think of this procedure as a proper algorithm, we must
be certain that it always terminates in a finite amount of time.
What would it mean for this computation to continue forever on a
given shape $S$?  If we succeed in building $n$-patches  as $n$
increases without bound, then we know that $S$ must tile the plane.
But because we never find a patch of copies of $S$ that tiles isohedrally,
we also know that $S$ cannot have a finite isohedral number, and therefore
cannot admit periodic tilings.  In other words, such an $S$
would have to be an \emph{aperiodic monotile}: a shape that admits
tilings, but never any that are periodic.

It is easy to construct a single tiling that is non-periodic; even
a humble $2\times 1$ brick admits uncountably many non-periodic
tilings.  Far more interesting is when the shape itself is
just flexible enough to admit tilings, but constrained enough not
to allow periodicity to arise.  This quality is so exotic that when
Hao Wang first proposed the notion of aperiodicity in the early
1960s, it was to conjecture immediately that sets of shapes that
behaved this way did not exist.  Wang's student Robert Berger later
found the first aperiodic set, comprising over 20,000 tiles, and
the race was on to achieve aperiodicity with fewer shapes.  Progress
was swift for about ten years, culminating in Penrose's discovery
of the kite and dart, a set of size two.  But despite a few close calls,
in the form of single shapes that tiled aperiodically under mild 
variations to the rules of the game, an aperiodic monotile could not
be found.

\begin{figure}
\begin{center}
\includegraphics[width=\textwidth]{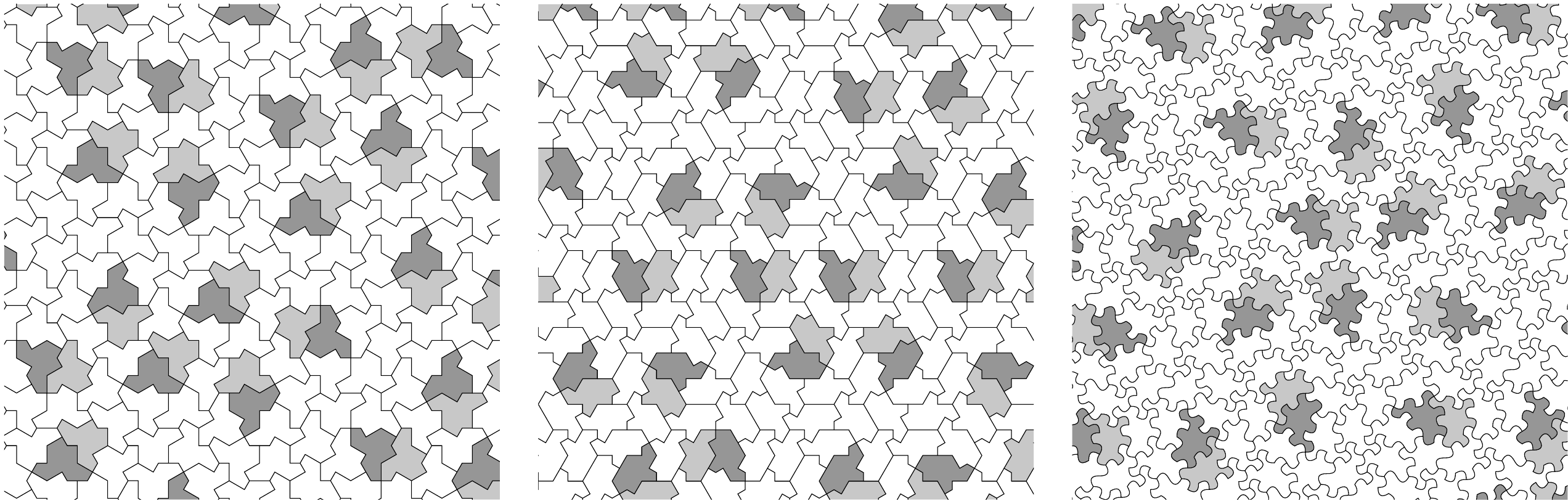}
\end{center}
\caption{\label{fig:aperiodic}Portions of tilings by 
	hats (left), turtles (centre), and a Spectre (right).}
\end{figure}

David Smith's discovery of the hat polykite in 2022 ended the
drought.  Fortunately, Chaim Goodman-Strauss, Joseph Myers and I
were well positioned to analyze a candidate aperiodic monotile, if
only one were to be proposed.  On top of our years of collective
study of the topic, Joseph and I had our software; we would have
been able to compute the hat's Heesch number if it were a non-tiler,
or its isohedral number if it tiled periodically.  I used my software
to compute a 16-patch of hats, confirming immediately that if the
hat did not tile the plane, then its Heesch number would have to
be at least 16, a preposterous advance beyond the record of six.
Assuming that the hat did admit tilings, David and I then scrutinized
my computer-generated patches to divine the structure of its tilings.
That process eventually yielded substitution rules, which can be
used to generate arbitrarily large patches of hats, showing that
the shape tiles the plane (Figure~\ref{fig:aperiodic}, left).  We
then generated exhaustive lists of small patches of hats that can
appear in infinite tilings.  We used those to rule out periodic
tilings by hats, completing a proof of the hat's aperiodicity~\cite{hat}.

In the course of our work, David found a second mysterious polykite
we call the ``turtle'', this one a union of ten kites (compared to
the hat's eight).  Joseph discovered that these two shapes were
connected by a continuum of shapes that all tiled aperiodically in
the same way (Figure~\ref{fig:aperiodic}, centre),
with the sole exception of an equilateral polygon
that we referred to prosaically as $\mathrm{Tile}(1,1)$.  He
then harnessed that continuum to drive a remarkable second proof
of the hat's aperiodicity that did not rely directly on computer
assistance.

When we first introduced the world to the hat, one objection was
raised more than any other.  Every tiling by hats must contain a
mix of left-handed and right-handed hats (i.e., a hat and its mirror
reflection).  Some people chose to regard these as two distinct
shapes, perhaps guided by the intuition of manufacturing physical
tiles.  While geometry in general, and tiling theory in particular,
support the point of view that left-handed and right-handed copies
of a shape are to be considered congruent, we were left with the
interesting question of whether aperiodicity could be achieved
without the use of reflections.  At first we thought that our work
could shed no light on that problem. However, David soon noticed
that miraculously, the previously discarded shape $\mathrm{Tile}(1,1)$
pointed the way to ``Spectres'' (Figure~\ref{fig:aperiodic}, right),
which are one-sided aperiodic monotiles: they admit only non-periodic
tilings in which all tiles have the same handedness~\cite{spectre}.
I remain dumbfounded by this second discovery.  As unlikely as the
hat was in the first place, I see no reason whatsoever why this
restricted form of aperiodicity should arise as a side effect of
that work.

While we now have a few shapes that prove that aperiodic monotiles
exist, we still have essentially no information about the broader
landscape of aperiodicity.  We do not know when and why aperiodicity
occurs, or where we might look for it.  Mathematics would not be
so cruel as to offer us just these few aperiodic monotiles, and so
the first task must be to find more of them.  Without a clear theory
to rely on, and with no further prophetic visions from the mind of
David Smith, we are left with no choice but to continue sifting
through collections of shapes like polyforms in search of a few
exquisite gems.  We should expect breakthroughs to be infrequent.
Joseph found no hints of aperiodicity among trillions of polyominoes,
polyhexes, and polyiamonds.  Apart from the sublime anomaly of the
hat and turtle, a subsequent search of 500 billion polykites uncovered
no other unusual behaviour.

We can imagine many other questions about the existence of aperiodic
monotiles under various modified or restricted conditions, besides
the prohibition of reflections.  It would be interesting to find
an aperiodic monotile with bileteral reflection symmetry, which
would render moot the question of handedness, or one with rotational
symmetry.  We could also consider restricted classes of shapes: we
might focus on polyominoes, for example, or more ambitiously find
a reason why aperiodicity cannot arise there.  It is also easy to
pose questions about the ``simplest'' possible aperiodic monotiles,
under various conceptions of simplicity. For example, the hat is a
13-sided polygon.  Can there exist an aperiodic monotile with fewer
sides?  The most we can currently say is that it would have to be
non-convex and have at least five sides.

Although most of this article concerns shapes in two dimensions, we
can consider these same questions in higher-dimensional spaces.  A
shape called the Schmitt-Conway-Danzer biprism is both an existence
proof and a cautionary tale about the definition of aperiodicity
in 3D space.  It tiles space without ever permitting translational
symmetries, but it admits tilings that contain a ``screw motion''
(a rotation about a axis composed with a translation parallel to
that axis), which can be repeated any number of times.  Some people
regard this screw motion as uncomfortably close to translation, and
demand a \emph{strongly aperiodic} 3D monotile, one that admits
tilings whose symmetries never include an infinite cyclic subgroup
of any kind.  Many of the ideas and algorithms we used to prove the
hat's aperiodicity could be adapted to work in 3D space, but
personally I am daunted by the prospect of deducing the behaviour
of any candidate shapes that might be discovered there.  In still
higher dimensions, it seems as if the phenomenon of aperiodicity becomes
more commonplace.  In recent work, Rachel Greenfeld and Terence Tao 
proved that in a sufficiently high number of dimensions---high enough
that the exact number is not known---there exist aperiodic monotiles
that tile by translation alone~\cite{GT2}.

Before we knew of the existence of aperiodic monotiles, I was drawn
to that problem much as I still am to Heesch numbers and isohedral
numbers.  Here again I savoured the equivocal tension that both
existence and non-existence of these shapes were real possibilities,
and hoped some day to learn which answer held.  I feel that we have
been granted a rare glimpse at the ``personality'' of the Euclidean
plane, as if it might have made up its mind either way on the
existence of aperiodic monotiles.  I'm glad it chose to permit them.

\section*{Final Thoughts}

The artist M.C. Escher observed that the construction of a tiling
requires a negotiation between every two tiles sharing an edge, as
they fight to divide up the space between them to satisfy artistic
(or mathematical) goals.  In tilings or patches by copies of a
single shape, that shape must somehow encode along its boundary all
the information it needs to determine the long-range behaviour of
the structures that it admits.  Viewed in that light, a shape like
Ba\v{s}i\'c's polygon with Heesch number six represents a marvellous
feat of geometric engineering: it is flexible enough to support six
layers of surrounding copies, but not seven.  Similarly, Joseph
Myers's polyhex with isohedral number 10 demands a certain minimum
investment of complexity in the assembly of a large patch, and then
settles down into a simple tiling by copies of that patch.  In both
cases, a single shape exerts a powerful influence over its widely
separated copies.  We have a great deal to learn about this particular
variety of ``spooky action at a distance''.

It is interesting to speculate about quantifying the amount of
information encoded in a shape's boundary, and the relationship
between that information and the kinds of patches and tilings the
shape admits.  One of the most remarkable things about the hat is
that its shape is so unassuming.  Somehow, though, the hat's boundary
contains just enough information to position it in the seemingly
infinitesimal gap between shapes that do not tile and shapes that
tile periodically. On the other hand, shapes with high Heesch numbers
and high isohedral numbers generally appear much more complex.
Perhaps we should regard this apparent complexity as evidence that
problems concerning Heesch numbers and isohedral numbers will prove
to be more ornery than finding an aperiodic monotile.  Or perhaps
we can draw the more profound conclusion that aperiodicity in some
forms will ultimately turn out to be a straightforward phenomenon, 
if only we can unlock a few more of its mysteries.

    
{\setlength{\baselineskip}{13pt} 
\raggedright				
\bibliographystyle{IEEEtran}
\bibliography{tilings}
} 

\end{document}